# Modeling macroeconomic time series via heavy tailed distributions

## J. A. D. Aston[1]

*Academia Sinica*

**Abstract:** It has been shown that some macroeconomic time series, especially those where outliers could be present, can be well modelled using heavy tailed distributions for the noise components. Methods for deciding when and where heavy-tailed models should be preferred are investigated. These investigations primarily focus on automatic methods for model identification and selection. Current methods are extended to incorporate a non-Gaussian selection element, and various different criteria for deciding on which overall model should be used are examined.

## 1. Introduction

While time series analysis is a rich topic for theoretical research, the implications of such work can impact many applied sciences, including physical, biological and social. An example of the application of statistical time series methods in economics is the seasonal adjustment of macroeconomic data, one of the primary functions carried out by many statistical agencies worldwide. These adjustments allow comparisons of economic indicators to be made in the presence of seasonal variations, and allow economic decisions to be made without the confounding factors of seasonal fluctuations. The seasonally adjusted data is an unobserved component in the data, and must be estimated using a model, which can be either parametrically or non-parametrically specified. However, the estimates from the model can be seriously affected by changes in the data, especially outliers in the data.

Methods have been developed to account for outliers in commonly used time series models for seasonal adjustment using heavier tailed distributions than the Gaussian [3, 5, 10], such as the t-distribution or mixtures of normals. The aim of this paper is to present findings regarding how to select whether a heavy tailed model is required for a data set based on the performance of several model selection criteria.

Currently, most statistical agencies use one of two packages for seasonal adjustment, X-12-ARIMA [6] from the US Census Bureau, or TRAMO/SEATS [7] from the Bank of Spain. These two programs both seasonally adjust the data, but in intrinsically different ways. X-12-ARIMA uses prespecified filters to remove the seasonal component from the data, in a non-parametric fashion, while TRAMO/SEATS uses the ARIMA methodology of [4] to determine the model and estimate the seasonal component. It will be this second approach that will be generalised to include non-Gaussian components. Recently, a new integrated version

[1]Institute of Statistical Science, Academia Sinica, 128 Academia Sinica, Sec 2, Taipei 11529, Taiwan, ROC, Tel: **+886-2-2783-5611 ext. 314** Fax: **+886-2-2783-1523** e-mail: `jaston@stat.sinica.edu.tw`

of the software has been released [11] and the methodology also applies to this package.

Two common types of outliers are additive outliers (shocks) and level shift (break) outliers. The first refers to a single data point that is out of character for the data given the model, whereas the second refers to a discontinuous jump either up or down in the level of the data at some point, and continuing on after that point for the rest of the series. These outliers lead to errors in the determination of the underlying components, and thus need to be accounted for. However, they are, by nature, not known prior to the modeling, and are functions of both the data and the model. A data point may be an outlier for one model but not another, and its status as an outlier may change when new data is added. This second feature is especially relevant for this study, as seasonal macroeconomic data is usually continuously updated month by month, or quarter by quarter.

Firstly, a very brief introduction to ARIMA models for seasonal adjustment will be given, followed by characterisations of the usual methods of outlier detection for both the Gaussian and non-Gaussian cases. While attempting to solve similar problems, the two approaches are appreciably different. Section 4 introduces the model selection criteria to be considered and also some justification for their usage. Section 5 provides some examples of real series where adjustment using heavy tailed models gives better performance than using Gaussian models and the final section provides discussion.

## 2. Seasonal adjustment and ARIMA model based decomposition

A seasonal time series $y_t$ can be expressed as the sum of unobserved components,

$$(1) \qquad y_t = S_t + N_t$$

where $S_t$ represents the seasonal component and $N_t$, the remaining non-seasonal component. Thus the seasonally adjusted series $y^{(sa)}$ is

$$(2) \qquad y_t^{(sa)} = y_t - S_t$$

and is also unobserved by the nature of being a function of $S_t$. ARIMA model based (AMB) decomposition specifies ARIMA models for the unobserved components and estimates these from both the data and from the overall model for the data as a whole.

Box and Jenkins [4] introduced a class of seasonal ARIMA models that model macroeconomic data well. The simplest model of this form is the airline model, which models differenced data as a product of moving average (MA) processes;

$$(3) \qquad (1-B)(1-B^s)y_t = (1-\theta B)(1-\Theta B^s)\epsilon_t$$

where $By_t = y_{t-1}$, $s$ represents the seasonal periodicity, and $\theta$ and $\Theta$ are the MA parameters associated with the non-seasonal and seasonal MA parts respectively. $\epsilon_t$ is assumed to be an iid Gaussian white noise process with variance $\sigma^2$. This model can be generalised by altering the degrees of the MA polynomials and adding Autoregressive (AR) parts to the left hand side of the equation. It will be assumed here that the differencing is not modified, and that it remains of an airline type. This includes the restriction that the MA and AR parameters cannot be unit roots as these would alter the overall differencing.

For simplicity, the airline model is considered explicitly, although other generalisations can be handled similarly. The AMB decomposition for the airline model, can be expressed in the following way using the decomposition of [9]

$$U(B)S_t = \theta_S(B)\omega_t$$
$$(1-B)^2 T_t = \theta_T(B)\eta_t$$
$$I_t = \varepsilon_t$$

(4)

where $U(B) = (1 + B + \cdots + B^{s-1})$ and $\omega_t, \eta_t, \epsilon_t$ are independent white noise processes and

$$y_t = S_t + T_t + I_t. \tag{5}$$

In order to define a unique solution (which may or may not exist), the restriction is taken that the pseudo-spectral densities of the seasonal and trend components have a minimum of zero (in line with the admissible decompositions of [9]). When this condition cannot be met without resulting in a negative variance for $I_t$, the decomposition is said to be inadmissible. Throughout the paper, only parameter combinations resulting in admissible decompositions will be assumed, which is almost always the case for macroeconomic data.

The parameter functions $\theta_T(), \theta_S()$ and the variances of $\omega_t, \eta_t, \epsilon_t$ are all functions of the underlying parameters $\theta, \Theta$ and $\sigma^2$. They can be calculated from the partial fraction decomposition of the pseudo-spectral densities, and the minimisation of each resulting component. This is usually done, for example in the SEATS software, after maximum likelihood estimation of the parameters has taken place, to give a final adjustment of the data.

## 2.1. Gaussian outlier adjustment

The TRAMO package [7] is the most widely used method for automatic model identification of seasonal ARIMA models for macroeconomic series. The program is used to estimate the order of differencing, the orders of the AR and MA components and also any outliers and common regressor effects that might be present. Inherently in this paper, it has been assumed that the order of differencing is the same as the airline model, but this assumption can be easily relaxed without significant change in the approaches outlined. All the other parts of the TRAMO procedure are used in exactly the same way in this paper as given in [8] except for the part relating to outlier detection.

The TRAMO software determines outliers as part of the automatic model identification portion of the program. Critical values for the thresholds at which data points are assumed to be outliers are chosen either by the user or from the length of the series. Outliers are found by determining whether the significance of the regression coefficients determined by assuming an outlier, be it an additive outlier or a level shift, has occurred at each point in the data. This is done iteratively, by adding in the largest regressor above the threshold (if one exists) and then repeating the exercise. A final check is made at the end to ensure that all regressors are still above the threshold for the final model.

This method effectively removes the data point when it is considered to be an outlier. When new data points arrive every month/quarter, the stability of the seasonal adjustment can heavily rely on the stability of the designated outliers to this new data as an outlying data point can then be added back in to the estimation if no longer classified as an outlier.

### 3. Heavy-tailed models to account for outliers in ARIMA component models

Aston and Koopman [3] proposed this alternate methodology to the Gaussian classification of outliers by using heavier tailed distributions to weight data points rather than making binary decisions on outliers. A short summary of the methodology is now given.

The decomposition model (4) can be modified to incorporate non-Gaussian components. In order to retain a similar structure, it is assumed that the components have the same variance as the decomposition would predict, but different densities are used to incorporate heavier tails.

The irregular component can be modified to include the t distribution in order to account for additive outliers as

$$(6) \qquad I_t^* \sim t(0, \sigma_I^2, \nu), \qquad t = 1, \ldots, n,$$

where $\nu > 2$ is the number of degrees of freedom and $\sigma_I^2$ is the variance, which is constant for any $\nu$. In the case of an irregular modelled by a mixture of normals,

$$(7) \qquad I_t^* \sim (1 - \rho)\mathcal{N}(0, \sigma_I^2) + \rho\mathcal{N}(0, \sigma_I^2 \lambda), \qquad t = 1, \ldots, n,$$

where $0 \leq \rho \leq 1$ determines the intensity of outliers in the series and $\lambda$ measures the magnitude of the outliers.

The decomposition model with a t-distributed irregular term can be expressed in its canonical form by

$$y_t = S_t + T_t + I_t^*, \qquad I_t^* \sim t(0, \sigma_I^2, \nu), \qquad t = 1, \ldots, n.$$

where $t(0, \sigma_I^2, \nu)$ refers to the t-density. This model has the same number of parameters as the original model specification except that the t density has one additional parameter (the degrees of freedom $\nu$) and the mixture of normals has two additional parameters (the intensity and the variance scalar).

To robustify the decomposition model against breaks in trend we consider the trend specification

$$(8) \qquad (1 - B)^2 T_t^* = \theta_T(B)\eta_t^*, \qquad \eta_t^* \sim t(0, \sigma_\eta^2, \nu_\eta),$$

where the t-distribution $t(0, \sigma_\eta^2, \nu_\eta)$ can be replaced by a mixture of normals distribution. The decomposition model with heavy tailed densities for both the trend innovations and the irregular is given by $y_t = S_t + T_t^* + I_t^*$ where the latter two components are given by (8) and (6), respectively.

These models can be estimated through the use of importance sampling as was described in [3]. For calculation, it is important to note that the decomposition must now be incorporated into the maximum likelihood estimation, as the individual components are modified, yet, they are still dependent on the overall ARIMA model for the series. However, fast algorithms for the decomposition make the estimation feasible.

### 4. Model selection

One of the most important issues is deciding which model to use and when. Here three different approaches are investigated, one based on the moments of the data, another using an empirical evaluation of a model selection criteria and the last based on the stability of the estimated components when the data is updated.

Although the model is complex and the estimation method of the maximum likelihood an approximation, certain properties of the model can be usefully investigated by considering a simpler model. Take the simple case of choosing between two noise models, one t-distributed and the other gaussian,

$$
\begin{aligned}
y_\tau = \epsilon_\tau \;\; \epsilon_\tau &\sim N(0, \phi_N) \\
y_\tau = \epsilon_\tau \;\; \epsilon_\tau &\sim t(0, \phi_T, \nu).
\end{aligned}
\tag{9}
$$

and all the $y_\tau$ are iid from one model or the other.

As can be seen, these are essentially nested models where the parameter of interest is the $\nu$ degrees of freedom parameter. Slightly different from the usual nesting setup, the gaussian model is the limiting distribution of the t-model as $\nu \to \infty$. A proof is given in the appendix to show by looking at a function of the moments of the data (essentially the kurtosis), a test can be performed as to whether the error term under investigation comes from a normal or t-distributed model. The test simply considers whether $\sqrt{n}(Z_n - 3)$ comes from $N(0, 24)$, as should be the case asymptotically, if the data are normally distributed, while for the t-distribution, the sequence will diverge. In the actual data case, the test will be applied to the irregular component data fitted under the normal model for parameter estimation.

As can be seen from the simulations in Figure 1, even for small samples of the size of the real data under investigation, there is a marked difference in the distribution of the statistic between the two models.

In addition to the model choice given above, two other methods are investigated. AIC [2] seems to be well suited to this problem, as the models are essentially nested. However, the theoretical justification for AIC requires that the maximum likelihood

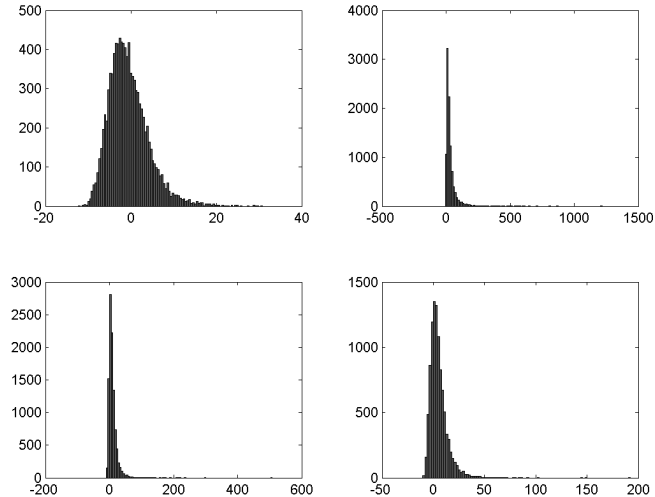

FIG 1. *Small sample kurtosis estimator distributions for four different models as generated from* $10^5$ *simulated samples of* $n = 150$. *(top left) Gaussian, (top right) t dist ($\nu = 5$), (bottom left) t dist ($\nu = 10$), (bottom right) t dist ($\nu = 15$)*

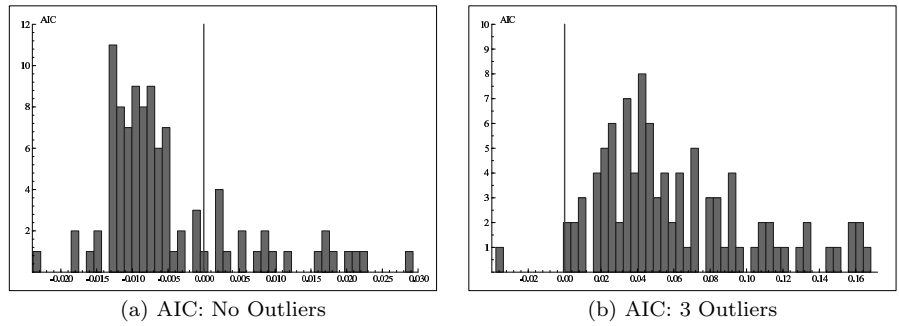

Fɪɢ 2. *AIC distributions for* 100 *series.* (a) *No outliers are present in the data* (b) *Three Outliers are present in the data. Black line corresponds to the value where AIC chooses Gaussian vs t-distributed model (left and right sides respectively).*

estimate of the parameters be an interior point within the parameter space, not on the boundary. Here, as the t distribution becomes the normal distribution as $\nu \to \infty$, this assumption is violated. However, in practice it can be seen that while the use of AIC might not yet be theoretically justified, in the case of the type of data under investigation, simulation results seem to be promising.

A small simulation study was carried out using the airline model (parameters of $\theta = 0.7$, $\Theta = 0.7$, $\sigma^2 = 1$) and two data sets generated, one where there were no outliers in the data, and one where there were three additive outliers added to the data (points shifted by 5 times the sd of the irregular component). Histograms of the AIC differences have been plotted in Figure 2. Corresponding histograms for AICc and BIC have also been generated (not shown) with similar results. AIC chooses the larger (t-distributed) model when there are outliers present in the data, and chooses the smaller model when outliers are not present, with an error rate of the same order as traditional AIC would predict.

In addition to the model selection criteria, an empirical measure was assessed for determining which model to use. This measure relies on using out-of-sample data to examine the estimates of the in-sample seasonally adjusted data when future data becomes available. A crude, yet seemingly promising, procedure is to withhold the final year of data, and to plot the changes in the seasonal components from the two samples, with and without the extra year of data. Given the problems of revisions when releasing macroeconomic data, adjusted series that remain stable when future data is added are to be preferred to adjusted series that change. By examining the plots of the differences, or some overall average change, such as the mean absolute difference between the two adjustments, the stability of the seasonal component to additional data can be quantified. Whilst this statistic is hard to justify theoretically given the complex nature of the model, it will be seen in the examples that it does seem to capture differences between the two approaches. Theoretical justification of this statistic will be the subject of future work.

## 5. Examples

Many macroeconomic series do not have large problems with outliers, and thus the methods described here will not be applicable. However, there are a sizeable proportion of series released by agencies such as the US Census Bureau where outliers do occur. When several series from the Census Bureau were investigated, two series where additive outliers seemed to be present were the Automobile Retail Series and the Material Handling Equipment Manufacturing Series (or u33mno as it is also known). These two series were analysed with the two different approaches, and the model selection criteria were used to determine which was the most appropriate model. Both series contained 155 data points (Feb 1992 until Dec 2004). For the seasonal stability plots, the data from 2004 was withheld in one analysis and used in the other and the results compared.

As can be seen in Figure 3, the seasonal stability plots for the automobile retail series indicate that the t-distribution provides a more stable adjustment than the normal model. This is chiefly because when an extra year of data is added, the number of outliers detected in the series changed. Thus this caused large changes in the seasonal pattern for the normal model. However, as there is no discrete detection process for the t-distributed model, there was a more continuous change in the seasonal pattern when the extra year of data was added, and thus the stability was greater. There was a low number of degrees of freedom (approximately 5-6) estimated for this model. The difference between the two models is also well detected by all the model selection criteria (Table 1). In addition, the moment estimator has a large value, well outside the 95% confidence interval range of a $N(0,24)$ and therefore normality of the error terms are rejected for both the 143 and 155 length series.

The same conclusion can be reached with Figure 4 for the u33mno series, although the seasonal pattern was more erratic for both models and thus the stability of the seasonal was closer in both models. This was also shown in smaller differences for the model selection criteria in Table 1, with all the criteria being borderline as to which model to use, especially given the small sample nature of

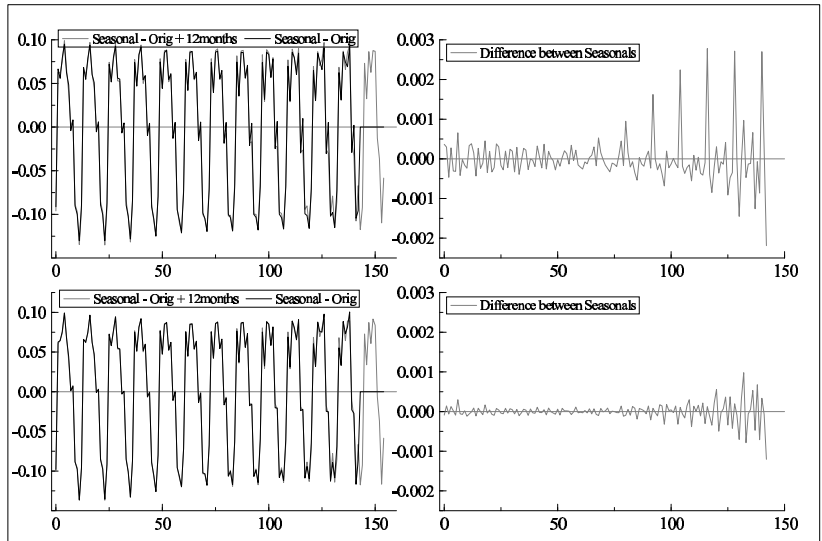

Seasonal Differences (Log Data). top row: Gaussian Data, bottom row: t-distributed model data, left column: Seasonal Patterns, right column: RMS Differences in Patterns when extra 12 observations added

FIG 3. *Automobile Retail Series from Feb 1992-Dec 2004 (US Census Bureau). This example shows that the seasonal difference plot finds a large change in the seasonal pattern with an extra year of data when a Gaussian model is used, but this change is reduced when using the model containing the t-distribution*

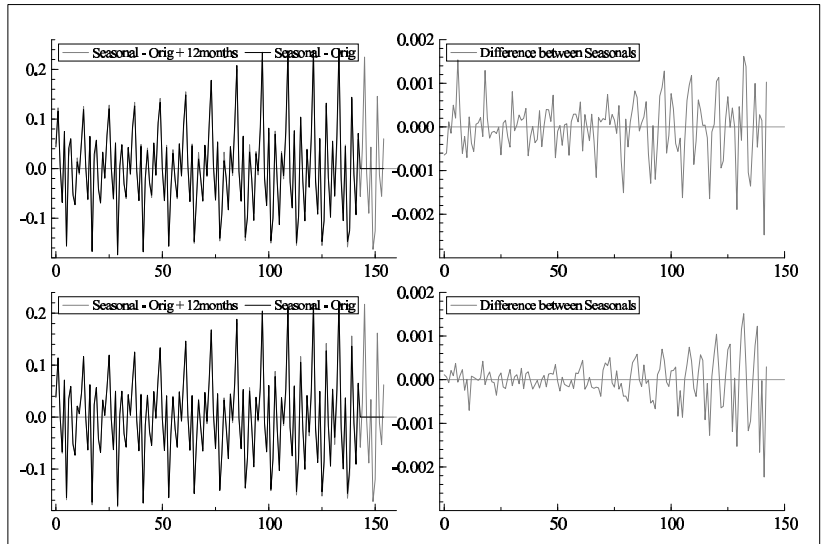

Seasonal Differences (Log Data). top row: Gaussian Data, bottom row: t-distributed model data, left column: Seasonal Patterns, right column: RMS Differences in Patterns when extra 12 observations added

Fig 4. *u33mno Series from Feb 1992-Dec 2004 (US Census Bureau). Again, this example shows that the seasonal difference plot finds a large change in the seasonal pattern with an extra year of data when a Gaussian model is used, but this change is reduced (slightly less than the auto retail series) when using the model containing the t-distribution.*

TABLE 1

*Model selection and comparison of the two example series. Bold face indicates the optimal value for selecting the model*

| Data | Model | Length | *df* | Sample Kurtosis | LogLik | AIC | AICc | BIC | Seas Mean Abs Diff |
|------|-------|--------|------|-----------------|--------|-----|------|-----|--------------------|
| Auto | G | 143 | - | 35.8 | 217.18 | -2.97 | -1.95 | -2.86 | - |
| Auto | T | 143 | 5.63 | - | 223.00 | **-3.04** | **-2.02** | **-2.91** | - |
| Auto | G | 155 | - | 47.5 | 236.75 | -2.99 | -1.97 | -2.89 | 0.0629 |
| Auto | T | 155 | 5.10 | - | 243.05 | **-3.06** | **-2.04** | **-2.94** | **0.0381** |
| u33mno | G | 143 | - | 11.1 | 92.29 | -1.22 | -0.20 | **-1.12** | - |
| u33mno | T | 143 | 9.17 | - | 93.59 | **-1.23** | **-0.21** | -1.10 | - |
| u33mno | G | 155 | - | 13.9 | 105.19 | -1.40 | -0.38 | **-1.30** | 0.0651 |
| u33mno | T | 155 | 7.43 | - | 107.14 | **-1.41** | **-0.39** | -1.29 | **0.0266** |

series. However, given the increase in the stability and the borderline nature, the t-distribution model will probably be preferred in the case of u33mno as well.

It can be noted in both Figures 3 and 4 that, for both the t-distributed and the Gaussian models, the instability within the estimates does increase towards the end of the series. This is due to estimates being weighted functions of other data. Data that is close to the point to be estimated (either directly or as a multiple of the seasonal period) is more heavily weighted than data that is further away. Thus when new data is added, the estimates towards the end of the series are more heavily affected than the estimates nearer the beginning of the series.

## 6. Discussion

In this paper, model selection criteria have been proposed to choose between Gaussian and heavier tailed distributions. Particular emphasis has been placed on choosing between the t-distribution and the Gaussian distribution for modeling the irregular component where additive outliers occur. However, all the techniques are generalisable to other distributions such as mixtures of normals and to the trend component to account for level shifts.

The model selection criteria have primarily been evaluated empirically for the data and models used in the paper. This is for two reasons. Firstly, the seasonal models under investigation are complex models, where the likelihood evaluation involves both approximations through importance sampling and also pseudo-spectral decomposition. Thus, results for these types of models are difficult to obtain explicitly. However, even for simpler models, only theoretical results have been obtained for the moment estimator selection procedure, given the nature of the model nesting, and the boundary problem. However, the results obtained from simulation are promising. This also suggests that theoretical justification of these and related results, which apply in many other modeling situations, may well be a worthwhile future research area.

### Acknowledgments

The author would like to express his gratitude to Ching-Kang Ing for all his help and suggestions with this work. He would also like to thank Siem Jan Koopman, Benedikt Pötscher and David Findley and an anonymous reviewer for their extremely helpful comments and discussions.

### Appendix A: Moment estimator

**Theorem 1** (Sample kurtosis). *Let $y_\tau$, $1 \leq \tau \leq n$ be a realisation from one of the two models given in* (9) *with finite positive variance.*

*Let*

$$(10) \qquad Z_n = \frac{\frac{1}{n}\sum(y_i - \bar{y})^4}{\left(\frac{1}{n}\sum(y_i - \bar{y})^2\right)^2}$$

*and if $\rightsquigarrow$ represents convergence in distribution then*

$$(11) \qquad \sqrt{n}(Z_n - 3) \rightsquigarrow N(0, 24)$$

*when $\nu \to \infty$ and diverges for $\nu$ finite.*

*Proof.* Both the denominator and numerator of $Z_n$ are moment estimators. If $y_\tau$ comes from the first model in (9) then following a similar method to [12, Example 3.5], let

$$\phi(a, b, c, d) = \frac{d - 4ca - 6ba^2 - 3a^4}{(b - a^2)^2}$$

then

$$Z_n = \phi(\overline{Y}, \overline{Y^2}, \overline{Y^3}, \overline{Y^4})$$

where $\overline{Y^j} = \frac{i=1}{n}\sum_1^n y_i^j$ and $\sqrt{n}(\overline{Y} - \alpha_1, \overline{Y^2} - \alpha_2, \overline{Y^3} - \alpha_3, \overline{Y^4} - \alpha_4)$ is asymptotically mean zero normal by the CLT where $\alpha_j$ is the $j$th moment of $y_1$ wlog.

If $X_i = \frac{y_i}{\phi_N}$, and using the fact that the odd moments of the standard normal are zero and the even moments given by $\frac{2n!}{2^n n!}$ (and thus the first eight moments are also finite),

$$\sqrt{n} \begin{pmatrix} \overline{X} \\ \overline{X^2} - 1 \\ \overline{X^3} \\ \overline{X^4} - 3 \end{pmatrix} \rightsquigarrow N \left( 0, \begin{pmatrix} 1 & 0 & 3 & 0 \\ 0 & 2 & 0 & 12 \\ 3 & 0 & 15 & 0 \\ 0 & 12 & 0 & 96 \end{pmatrix} \right).$$

The function $\phi$ is differentiable at the point $(\alpha_1, \alpha_2, \alpha_3, \alpha_4) = (0, 1, 0, 3)$, and equals $(0, -6, 0, 1)$. Hence, by use of the delta method

$$\sqrt{n}(Z_n - 3) \rightsquigarrow N(0, 24)$$

If $y_\tau$ comes from the second model in (9) (and assuming $\nu > 4$) then

$$\frac{1}{n} \sum (y_i - \bar{y})^4 \to 3\phi_T^2 \frac{\nu^2}{(\nu - 2)(\nu - 4)} \quad (a.s.)$$

and

$$\left( \frac{1}{n} \sum (y_i - \bar{y})^2 \right) \to \phi_T \frac{\nu}{\nu - 2} \quad (a.s.)$$

by explicit calculation of the moments of the t-distribution [1] and thus

$$Z_n \to 3. \frac{(\nu - 2)}{(\nu - 4)} \quad (a.s.).$$

As $y_\tau$ comes from second model of (9), $\nu$ is finite and as $n \to \infty$

$$\sqrt{n}(Z_n - 3)$$

will diverge. $\qquad\qquad\qquad\qquad\qquad\qquad\qquad\qquad\qquad\qquad\qquad\qquad\qquad\qquad\Box$